\documentclass[10pt]{amsart}
\usepackage{amscd}
\usepackage{amssymb}
\usepackage[all]{xy}
\title{Erratum to: Deformation Quantization in Algebraic Geometry}
\author{Amnon Yekutieli}
\address{A. Yekutieli: Department of  Mathematics
Ben Gurion University, Be'er Sheva 84105, Israel}
\email{amyekut@math.bgu.ac.il}
\date{13 August 2007}

\newtheorem{thm}[equation]{Theorem}
\newtheorem{cor}[equation]{Corollary}

\newtheorem{lem}[equation]{Lemma}
\theoremstyle{definition}

\newtheorem{exa}[equation]{Example}

\numberwithin{equation}{section}

\newcommand{\iso}{\xrightarrow{\simeq}}

\newcommand{\opn}{\operatorname}
\newcommand{\cat}[1]{\operatorname{\mathsf{#1}}}

\newcommand{\rmitem}[1]{\item[\text{\textup{(#1)}}]}
\newcommand{\mfrak}[1]{\mathfrak{#1}}
\newcommand{\mcal}[1]{\mathcal{#1}}
\newcommand{\msf}[1]{\mathsf{#1}}
\newcommand{\mbf}[1]{\mathbf{#1}}
\newcommand{\mrm}[1]{\mathrm{#1}}
\newcommand{\mbb}[1]{\mathbb{#1}}

\newcommand{\smfrac}[2]{{\textstyle \frac{#1}{#2}}}
\newcommand{\tup}[1]{\textup{#1}}
\newcommand{\bsym}[1]{\boldsymbol{#1}}
\newcommand{\boplus}{\bigoplus\nolimits}

\newcommand{\what}[1]{\widehat{#1}}

\newcommand{\set}[1]{\{ #1 \}}

\newcommand{\K}{\mbb{K} \hspace{0.05em}}
\renewcommand{\d}{\mrm{d}}

\newcommand{\sprod}{{\textstyle \prod}}
\newcommand{\hatotimes}[1]{\, \what{\otimes}_{#1} \,}

\begin{document}

\maketitle

\begin{abstract}
This note contains a correction of the {\em proofs} of the main
results of the paper 
[A.\ Yekutieli, Deformation quantization in algebraic geometry,
Adv.\ Math.\ {\bf 198} (2005), 383-432]. The results are correct as
originally stated. 
\end{abstract}

\setcounter{section}{-1}
\section{Introduction}
 
This note contains a correction of the {\em proofs} of the main
results of \cite{Ye1}, namely Theorems 0.1 and 0.2. The results are correct as 
originally stated. 

The mistake in my original proofs was discovered Michel Van den
Bergh, and I thank him for calling my attention to it. The way to fix the
proofs is essentially contained in his paper \cite{VdB}.

Let me begin by explaining the mistake. As can be seen in Example 
\ref{exa1} below, the mistake itself is of a rather elementary nature, but it
was obscured by the  complicated context. 

Suppose $\K$ is a field of characteristic $0$, and $X$ is a smooth separated
\linebreak $n$-dimensional scheme over $\K$. Recall that the coordinate bundle 
$\opn{Coor} X$ is an infinite dimensional bundle over $X$, with free action by
the group $\mrm{GL}_{n, \K}$. 
The quotient bundle is by definition 
\[ \opn{LCC} X := \opn{Coor} X / \mrm{GL}_{n, \K} , \]
and the projection
$\pi_{\mrm{gl}} : \opn{Coor} X \to \opn{LCC} X$
is a $\mrm{GL}_{n, \K}$-torsor.

The erroneous (implicit) assertion in \cite{Ye1} is that the de Rham
complexes satisfy
\[ (\pi_{\mrm{gl} *}\, \Omega_{\opn{Coor} X})^{\mrm{GL}_{n}(\K)} = 
\Omega_{\opn{LCC} X} . \]
 From that I deduced (incorrectly, top of page 424) that the
Maurer-Cartan form 
$\omega_{\mrm{MC}}$ is a global section of the sheaf
\[ \Omega^1_{\opn{LCC} X} 
\, \what{\otimes}_{\mcal{O}_{\opn{LCC} X}} \,
\pi_{\mrm{lcc}}^{\what{*}}(\mcal{P}_X \otimes_{\mcal{O}_{X}}
\mcal{T}^{0}_{\mrm{poly}, X})  . \]
(This false, as can be seen from \cite[Lemma 6.5.1]{VdB}). 
This led to many incorrect formulas in \cite[Section 7]{Ye1}. 

The correct thing to do is to work with the infinitesimal action of the Lie
algebra
$\mfrak{g} := \mfrak{gl}_n(\K)$.
For $v \in \mfrak{g}$ one has
the contraction (inner derivative) $\iota_v$, 
which is a degree $-1$ derivation of the de Rham complex
$\pi_{\mrm{gl} *}\, \Omega_{\opn{Coor} X}$.
Recall that the Lie derivative is
$\mrm{L}_v := \d \circ \iota_v + \iota_v \circ \d$. 
A local section $\omega \in \pi_{\mrm{gl} *}\, \Omega_{\opn{Coor} X}$ 
is said to be $\mfrak{g}$-invariant if
$\iota_v(\omega) = \mrm{L}_v(\omega) = 0$
for all $v \in \mfrak{g}$. According to \cite[Lemma 9.2.3]{VdB} one has
\[ (\pi_{\mrm{gl} *}\, \Omega_{\opn{Coor} X})^{\mfrak{g}} = 
\Omega_{\opn{LCC} X} . \]

It is worthwhile to note that in my incorrect proof there was no need to invoke
Kontsevich's property (P5) from \cite{Ko}. The correct proof does require
property (P5) -- cf.\ \cite[Lemma 9.2.1]{VdB}.

\begin{exa} \label{exa1}
Here is a simplified example. Suppose $G$ is the affine algebraic group
$\mrm{GL}_{1, \K} = \opn{Spec} \K[t, t^{-1}]$, and 
$X$ is the variety $G$, with regular left action. The group of rational
points is $G(\K) = \K^{\times}$.
The action of $G$ on $X$ is free, the invariant ring is
$\mcal{O}(X)^{G(\K)} = \K$, and the quotient is 
$X / G = \opn{Spec} \K$.
For the de Rham complex 
\[ \Omega(X) = \mcal{O}(X) \oplus \Omega^1(X) = 
\K[t, t^{-1}] \oplus \K[t, t^{-1}] \cdot \d t  \]
we have
$\Omega(X)^{G(\K)} \neq \K$,
since it contains $t^{-1} \d t$. 
But for the infinitesimal action of the Lie algebra 
$\mfrak{g} := \mfrak{gl}_1(\K)$ it is easy to see that 
$\Omega(X)^{\mfrak{g}} = \K$.
\end{exa}

After some deliberation I decided that the best way to present the
erratum is by completely rewriting \cite[Section 7]{Ye1}. This is Section 1
below. Section 2 contains some additional minor corrections to
\cite{Ye1}.

\section{The Global $\mrm{L}_{\infty}$ Quasi-isomorphism} 

This is a revised version of \cite[Section 7]{Ye1}.
In this section we prove the main results of the paper \cite{Ye1}, namely
Theorem 0.1 (which is repeated here as Corollary \ref{cor9.1}),
and Theorem 0.2 (which is repeated here, with more details, as Theorem
\ref{thm5.1}). Throughout $\K$ is a field
containing $\mbb{R}$, and $X$ is a smooth
irreducible separated $n$-dimensional scheme over $\K$. 
We use all notation, definitions and results of \cite[Sections 1-6]{Ye1}
freely.  However the bibliography references relate to the list at
the end of this note.

Suppose $\bsym{U} = \{ U_{0}, \ldots, U_{m} \}$ is an open covering 
of the scheme $X$, consisting of affine open sets, each admitting 
an \'etale coordinate system, namely an \'etale morphism 
$U_{i} \to \mbf{A}^n_{\K}$. For every $i$ let
$\sigma_{i} : U_{i} \to \opn{LCC} X$ be the corresponding 
section of $\pi_{\mrm{lcc}} : \opn{LCC} X \to X$, 
and let $\bsym{\sigma}$ be the 
resulting simplicial section (see \cite[Theorem 6.5]{Ye1}).

Let $\mcal{M}$ be a bounded below complex of 
quasi-coherent $\mcal{O}_{X}$-modules.
The mixed resolution $\opn{Mix}^{}_{\bsym{U}}(\mcal{M})$ 
was defined in \cite[Section 6]{Ye1}. For any integer $i$
let
\[ \mrm{G}^i \opn{Mix}^{}_{\bsym{U}}(\mcal{M}) :=
\boplus_{j = i}^{\infty} \opn{Mix}^{j}_{\bsym{U}}(\mcal{M}) , \]
so 
$\{ \mrm{G}^i \opn{Mix}^{}_{\bsym{U}}(\mcal{M}) \}_{i \in \mbb{Z}}$
is a descending filtration of 
$\opn{Mix}^{}_{\bsym{U}}(\mcal{M})$
by subcomplexes, with 
$\mrm{G}^i \opn{Mix}^{}_{\bsym{U}}(\mcal{M}) = 
\opn{Mix}^{}_{\bsym{U}}(\mcal{M})$
for $i \leq 0$, and
$\bigcap_{i}\, \mrm{G}^i \opn{Mix}^{}_{\bsym{U}}(\mcal{M}) = 0$.
Let 
\[ \opn{gr}_{\mrm{G}}^i \opn{Mix}^{}_{\bsym{U}}(\mcal{M}) :=
\mrm{G}^i \opn{Mix}^{}_{\bsym{U}}(\mcal{M}) \ / \
\mrm{G}^{i+1} \opn{Mix}^{}_{\bsym{U}}(\mcal{M}) \]
and
$\opn{gr}_{\mrm{G}} \opn{Mix}^{}_{\bsym{U}}(\mcal{M}) :=
\boplus_{i}\, \opn{gr}_{\mrm{G}}^i 
\opn{Mix}^{}_{\bsym{U}}(\mcal{M})$.

By \cite[Proposition 6.3]{Ye1}, if $\mcal{G}_X$ is either 
$\mcal{T}^{}_{\mrm{poly}, X}$ or 
$\mcal{D}^{}_{\mrm{poly}, X}$, then 
$\opn{Mix}^{}_{\bsym{U}}(\mcal{G}_X)$ is a sheaf of DG Lie 
algebras on $X$, and the inclusion
\[ \eta_{\mcal{G}} : \mcal{G}_X \to \opn{Mix}^{}_{\bsym{U}}(\mcal{G}_X) \]
is a DG Lie algebra quasi-isomorphism. 

Note that if 
$\phi : \opn{Mix}^{}_{\bsym{U}}(\mcal{M}) \to
\opn{Mix}^{}_{\bsym{U}}(\mcal{N})$
is a homomorphism of complexes that respects the filtration
$\{ \mrm{G}^i \opn{Mix}^{}_{\bsym{U}} \}$,
then there exists an induced homomorphism of complexes
\[ \opn{gr}_{\mrm{G}}(\phi) : 
\opn{gr}_{\mrm{G}} \opn{Mix}^{}_{\bsym{U}}(\mcal{M}) \to
\opn{gr}_{\mrm{G}} \opn{Mix}^{}_{\bsym{U}}(\mcal{N}) . \]

Suppose $\mcal{G}$ and $\mcal{H}$ are sheaves of DG Lie algebras 
on a topological space $Y$. An $\mrm{L}_{\infty}$ morphism
$\Psi : \mcal{G} \to \mcal{H}$ is a sequence of sheaf morphisms
$\psi_j : \sprod^j \mcal{G} \to \mcal{H}$, such that for every 
open set $V \subset Y$ the sequence 
$\{ \Gamma(V, \psi_j) \}_{j \geq 1}$ is an $\mrm{L}_{\infty}$ 
morphism
$\Gamma(V, \mcal{G}) \to \Gamma(V, \mcal{H})$. If 
$\psi_1 : \mcal{G} \to \mcal{H}$ is a 
quasi-isomorphism then $\Psi$ is called an $\mrm{L}_{\infty}$ 
quasi-isomorphism.

Recall that there is a canonical quasi-isomorphism of complexes of
$\mcal{O}_X$-modules
\begin{equation} \label{eqn13.2}
\mcal{U}_1 : \mcal{T}_{\mrm{poly}, X} \to \mcal{D}_{\mrm{poly}, X} .
\end{equation}
According to \cite[Theorem 4.17]{Ye2}, the induced homomorphism 
\[ \opn{gr}_{\mrm{G}}(\opn{Mix}^{}_{\bsym{U}}(\mcal{U}_1)) :
\opn{gr}_{\mrm{G}} \opn{Mix}^{}_{\bsym{U}}(\mcal{T}_{\mrm{poly}, X}) 
\to \opn{gr}_{\mrm{G}} \opn{Mix}^{}_{\bsym{U}}(\mcal{D}_{\mrm{poly}, X}) \]
is a quasi-isomorphism.

\begin{thm} \label{thm5.1}
Let $X$ be an irreducible smooth separated $\K$-scheme.
Let $\bsym{U} = \{ U_{0}, \ldots, U_{m} \}$ 
be an open covering of $X$ consisting of affine open sets, each admitting 
an \'etale coordinate system, and let
$\bsym{\sigma}$ be the associated simplicial section of the bundle
$\opn{LCC} X \to X$. Then there is an 
induced $\mrm{L}_{\infty}$ quasi-isomorphism
\[ \Psi_{\bsym{\sigma}} = 
\set{\Psi_{\bsym{\sigma}; j}}_{j \geq 1} :
\opn{Mix}^{}_{\bsym{U}}(\mcal{T}_{\mrm{poly}, X}) 
\to \opn{Mix}^{}_{\bsym{U}}(\mcal{D}_{\mrm{poly}, X}) . \]
The homomorphism $\Psi_{\bsym{\sigma}; 1}$ respects the filtration
$\{ \mrm{G}^i \opn{Mix}^{}_{\bsym{U}} \}$,
and 
\[ \opn{gr}_{\mrm{G}}(\Psi_{\bsym{\sigma}; 1}) = 
\opn{gr}_{\mrm{G}}(\opn{Mix}^{}_{\bsym{U}}(\mcal{U}_{1})) . \]
\end{thm}

\begin{proof}
Let $Y$ be some
$\K$-scheme, and denote by $\K_Y$ the constant sheaf. For any $p$ 
we view $\Omega_Y^p$ as a discrete inv $\K_Y$-module, and we put 
on $\Omega_Y = \boplus_{p \in \mbb{N}}\, \Omega_Y^p$ 
direct sum dir-inv structure. So $\Omega_Y$ is a discrete (and 
hence complete)
DG algebra in $\cat{Dir} \cat{Inv} \cat{Mod} \K_Y$.

We shall abbreviate
$\mcal{A} := \Omega_{\opn{Coor} X}$, so that 
$\mcal{A}^0 = \mcal{O}_{\opn{Coor} X}$ etc.
As explained above, $\mcal{A}$ is a DG algebra in 
$\cat{Dir} \cat{Inv} \cat{Mod} \K_{\opn{Coor} X}$,
with discrete (but not trivial) dir-inv module structure.

There are sheaves of DG Lie algebras 
$\mcal{A} \, \what{\otimes} \
\mcal{T}^{}_{\mrm{poly}}(\K[[\bsym{t}]])$
and 
$\mcal{A} \, \what{\otimes} \
\mcal{D}^{}_{\mrm{poly}}(\K[[\bsym{t}]])$
on the scheme $\opn{Coor} X$. The differentials are
$\mrm{d}_{\mrm{for}} = \mrm{d} \otimes \bsym{1}$
and $\mrm{d}_{\mrm{for}} + \bsym{1} \otimes \mrm{d}_{\mcal{D}}$
respectively. As explained just prior to \cite[Theorem 3.16]{Ye1}, 
$\mcal{U}$ extends to a continuous $\mcal{A}$-multilinear 
$\mrm{L}_{\infty}$ morphism
\[ \mcal{U}_{\mcal{A}} = 
\{ \mcal{U}_{\mcal{A} ;j} \}_{j \geq 1} :
\mcal{A} \, \what{\otimes} \
\mcal{T}^{}_{\mrm{poly}}(\K[[\bsym{t}]]) \to
\mcal{A} \, \what{\otimes} \
\mcal{D}^{}_{\mrm{poly}}(\K[[\bsym{t}]])  \]
of sheaves of DG Lie algebras on $\opn{Coor} X$.

The MC form $\omega := \omega_{\mrm{MC}}$ is a global section of
$\mcal{A}^1 \, \what{\otimes} \
\mcal{T}^{0}_{\mrm{poly}}(\K[[\bsym{t}]])$
satisfying the MC equation in the DG Lie algebra
$\mcal{A} \, \what{\otimes} \
\mcal{T}^{}_{\mrm{poly}}(\K[[\bsym{t}]])$. 
See \cite[Proposition 5.9]{Ye1}.
According to  \cite[Theorem 3.16]{Ye1}, the global section 
$\omega' := \mcal{U}_{\mcal{A} ;1}(\omega) \in 
\mcal{A}^1 \hatotimes{} 
\mcal{D}^{0}_{\mrm{poly}}(\K[[\bsym{t}]])$
is a solution of the MC equation in the DG Lie algebra \linebreak
$\mcal{A} \, \what{\otimes} \
\mcal{D}^{}_{\mrm{poly}}(\K[[\bsym{t}]])$,
and there is a continuous $\mcal{A}$-multilinear
$\mrm{L}_{\infty}$ morphism
\[ \mcal{U}_{\mcal{A}, \omega} = \{ \mcal{U}_{\mcal{A}, \omega ;j} 
\}_{j \geq 1} : 
\big( \mcal{A} \, \what{\otimes} \,
\mcal{T}^{}_{\mrm{poly}}(\K[[\bsym{t}]]) \big)_{\omega} \to 
\big( \mcal{A} \, \what{\otimes} \,
\mcal{D}^{}_{\mrm{poly}}(\K[[\bsym{t}]]) \big)_{\omega'} \]
between the twisted DG Lie algebras. 
The formula is
\begin{equation} \label{eqn9.6}
\mcal{U}_{\mcal{A}, \omega; j}(\gamma_1 \cdots \gamma_j) = 
\sum_{k \geq 0} \smfrac{1}{(j+k)!}
\mcal{U}_{\mcal{A} ;j + k} (\omega^k \cdot \gamma_1 \cdots \gamma_j)
\end{equation} 
for 
$\gamma_1, \ldots, \gamma_j
\in \mcal{A} \, \what{\otimes} \, 
\mcal{T}^{}_{\mrm{poly}}(\K[[\bsym{t}]])$.
The two twisted DG Lie algebras have differentials
$\mrm{d}_{\mrm{for}} + \opn{ad}(\omega)$
and
$\mrm{d}_{\mrm{for}} + \opn{ad}(\omega')
+ \bsym{1} \otimes \mrm{d}_{\mcal{D}}$
respectively. 

This sum in (\ref{eqn9.6}) is actually finite, the number of nonzero terms in
it depending on the bidegree of $\gamma_1 \cdots \gamma_j$. 
Indeed, if
$\gamma_1 \cdots \gamma_j \in
\mcal{A}^{q} \, \what{\otimes} \, 
\mcal{T}^{p}_{\mrm{poly}}(\K[[\bsym{t}]])$,
then
\begin{equation} \label{eqn11.1}
\mcal{U}_{\mcal{A}; j+k} ( \omega^k \cdot \gamma_1 \cdots \gamma_j) \in 
\mcal{A}^{q+k} \, \what{\otimes} \, 
\mcal{D}^{p+1-j-k}_{\mrm{poly}}(\K[[\bsym{t}]]) ,
\end{equation}
which is is zero for $k > p-j+2$; see proof of  
\cite[Theorem 3.23]{Ye2}.

By \cite[Theorem 5.6]{Ye1} (the universal Taylor expansions) 
there are canonical isomorphisms of graded Lie algebras
in $\cat{Dir} \cat{Inv} \cat{Mod} \K_{\opn{Coor} X}$
\begin{equation} \label{eqn12.4}
\mcal{A} \, \what{\otimes} \,
\mcal{T}^{}_{\mrm{poly}}(\K[[\bsym{t}]]) \cong
\mcal{A}
\, \what{\otimes}_{\mcal{A}^0} \,
\pi_{\mrm{coor}}^{\what{*}}(\mcal{P}_X \otimes_{\mcal{O}_{X}}
\mcal{T}^{}_{\mrm{poly}, X}) 
\end{equation}
and
\begin{equation} \label{eqn12.5}
\mcal{A} \, \what{\otimes} \,
\mcal{D}^{}_{\mrm{poly}}(\K[[\bsym{t}]]) \cong
\mcal{A} \, \what{\otimes}_{\mcal{A}^0} \,
\pi_{\mrm{coor}}^{\what{*}}(\mcal{P}_X \otimes_{\mcal{O}_{X}}
\mcal{D}^{}_{\mrm{poly}, X}) . 
\end{equation}
\cite[Proposition 5.8]{Ye1} tells us that 
\[ \mrm{d}_{\mrm{for}} + 
\opn{ad}(\omega) = \nabla_{\mcal{P}}  \]
under these identifications. Therefore 
\begin{equation} \label{eqn13.1}
\begin{aligned}
& \mcal{U}_{\mcal{A}, \omega} =
\{ \mcal{U}_{\mcal{A}, \omega ;j} \}_{j \geq 1}  : 
\mcal{A} \, \what{\otimes}_{\mcal{A}^0} \,
\pi_{\mrm{coor}}^{\what{*}}(\mcal{P}_X \otimes_{\mcal{O}_{X}}
\mcal{T}^{}_{\mrm{poly}, X}) \\
& \hspace{10ex} \to \mcal{A} \, \what{\otimes}_{\mcal{A}^0} \,
\pi_{\mrm{coor}}^{\what{*}}(\mcal{P}_X \otimes_{\mcal{O}_{X}}
\mcal{D}^{}_{\mrm{poly}, X})
\end{aligned}
\end{equation}
is a continuous $\mcal{A}$-multilinear
$\mrm{L}_{\infty}$ morphism between these DG Lie algebras, 
whose differentials are
$\nabla_{\mcal{P}}$ and
$\nabla_{\mcal{P}} + \bsym{1} \otimes \d_{\mcal{D}}$ 
respectively.
As in the proof of \cite[Theorem 5.6]{Ye1}, under the identifications
(\ref{eqn12.4}) and (\ref{eqn12.5})  we have the equality
\begin{equation} \label{eqn12.6}
\mcal{U}_{\mcal{A}; 1} = 
\bsym{1} \otimes \pi_{\mrm{coor}}^{\what{*}}(\bsym{1} \otimes
\mcal{U}_1) ,
\end{equation}
i.e.\ it is the pullback of the map (\ref{eqn13.2}).

Let us filter the DG algebra $\mcal{A}$ by the descending filtration
$\{ \mrm{G}^i \mcal{A} \}_{i \in \mbb{Z}}$, where
$\mrm{G}^i \mcal{A} :=
\boplus_{j = i}^{\infty} \mcal{A}^i$. 
The DG Lie algebras appearing in equation (\ref{eqn13.1}) inherit
this filtration. From formulas (\ref{eqn9.6}) and (\ref{eqn11.1}) we see that
the homomorphism of complexes 
$\mcal{U}_{\mcal{A}, \omega; 1}$ respects the filtration, and from
(\ref{eqn12.6}) we see that 
\[ \opn{gr}_{\mrm{G}}(\mcal{U}_{\mcal{A}, \omega; 1}) = 
\opn{gr}_{\mrm{G}}(\mcal{U}_{\mcal{A}; 1}) =
\bsym{1} \otimes \pi_{\mrm{coor}}^{\what{*}}(\bsym{1} \otimes
\mcal{U}_1) . \]

Let $n := \opn{dim} X$. 
As noted earlier, the action of 
$\mfrak{g} := \mfrak{gl}_n(\K)$ gives
\[ (\pi_{\mrm{gl} *}\, \mcal{A})^{\mfrak{g}} = 
(\pi_{\mrm{gl} *}\, \Omega_{\opn{Coor} X})^{\mfrak{g}} = 
\Omega_{\opn{LCC}X} . \]
According to \cite[Lemma 9.2.1]{VdB}, the $\mrm{L}_{\infty}$ morphism 
$\mcal{U}_{\mcal{A}, \omega}$
commutes with the action of the Lie algebra $\mfrak{g}$. 
Therefore 
$\mcal{U}_{\mcal{A}, \omega}$ descends (i.e.\ restricts)
to a continuous
$\Omega_{\opn{LCC}X}$-multilinear $\mrm{L}_{\infty}$ morphism
\begin{equation} \label{eqn12.1}
\begin{aligned}
& \mcal{U}_{\mcal{A}, \omega}^{\mfrak{g}} : 
\Omega^{}_{\opn{LCC} X}
\, \what{\otimes}_{\mcal{O}_{\opn{LCC} X}} \,
\pi_{\mrm{lcc}}^{\what{*}}(\mcal{P}_X \otimes_{\mcal{O}_{X}}
\mcal{T}^{}_{\mrm{poly}, X}) \\
& \qquad \qquad \qquad \qquad \to \
\Omega^{}_{\opn{LCC} X}
\, \what{\otimes}_{\mcal{O}_{\opn{LCC} X}} \,
\pi_{\mrm{lcc}}^{\what{*}}(\mcal{P}_X \otimes_{\mcal{O}_{X}}
\mcal{D}^{}_{\mrm{poly}, X}) .
\end{aligned}
\end{equation}

The DG Lie algebras in formula (\ref{eqn12.1}) also have  
filtrations $\{ \mrm{G}^j \}_{j \in \mbb{Z}}$, the homomorphism 
$\mcal{U}_{\mcal{A}, \omega; 1}^{\mfrak{g}}$ respects this filtration, and we
now have
\begin{equation} \label{eqn12.7}
\opn{gr}_{\mrm{G}}(\mcal{U}_{\mcal{A}, \omega; 1}^{\mfrak{g}}) = 
\opn{gr}_{\mrm{G}}(\mcal{U}_{\mcal{A}; 1}^{\mfrak{g}}) =
\bsym{1} \otimes \pi_{\mrm{lcc}}^{\what{*}}(\bsym{1} \otimes
\mcal{U}_1) .
\end{equation}

According to \cite[Theorem 6.4]{Ye1}  there are induced operators 
\[ \Psi_{\bsym{\sigma}; j} := 
\bsym{\sigma}^*(\mcal{U}_{\mcal{A}, \omega; j}^{\mfrak{g}}) :
\prod\nolimits^j
\opn{Mix}^{}_{\bsym{U}}
(\mcal{T}^{}_{\mrm{poly}, X}) 
\to \opn{Mix}^{}_{\bsym{U}}
(\mcal{D}^{}_{\mrm{poly}, X}) \]
for $j \geq 1$. The $\mrm{L}_{\infty}$ identities in 
\cite[Definition 3.7]{Ye1},
when applied to the $\mrm{L}_{\infty}$ morphism 
$\mcal{U}_{\mcal{A}, \omega}^{\mfrak{g}}$, are of the form considered in 
\cite[Theorem 6.4(iii)]{Ye1}. 
Therefore these identities are preserved by $\bsym{\sigma}^*$,
and we conclude that the sequence 
$\Psi_{\bsym{\sigma}}  = 
\{ \Psi_{\bsym{\sigma}, j} \}_{j = 1}^{\infty}$
is an $\mrm{L}_{\infty}$ morphism.
Furthermore, $\Psi_{\bsym{\sigma}; 1}$ respects the filtration
$\{ \mrm{G}^i \opn{Mix}^{}_{\bsym{U}} \}$, and from (\ref{eqn12.7}) we get
\begin{equation} \label{eqn12.8}
\opn{gr}_{\mrm{G}}(\Psi_{\bsym{\sigma}; 1}) = 
\opn{gr}_{\mrm{G}}(\bsym{\sigma}^*(\mcal{U}_{\mcal{A}; 1}^{\mfrak{g}})) =
\opn{gr}_{\mrm{G}}( \opn{Mix}^{}_{\bsym{U}}(\mcal{U}_1)) . 
\end{equation}

According to \cite[Theorem 4.17]{Ye2}
the homomorphism 
$\opn{gr}_{\mrm{G}}(\opn{Mix}^{}_{\bsym{U}}(\mcal{U}_1))$ is a 
quasi-isomorphism. Since the complexes 
$\opn{Mix}^{}_{\bsym{U}}(\mcal{T}^{}_{\mrm{poly}, X})$
and
$\opn{Mix}^{}_{\bsym{U}}(\mcal{D}^{}_{\mrm{poly}, X})$
are bounded below, and the filtration is nonnegative and exhaustive,
it follows that $\Psi_{\bsym{\sigma}; 1}$ 
is also a quasi-isomorphism.
\end{proof}

\begin{cor} \label{cor11.1}
Taking global sections in Theorem \tup{\ref{thm5.1}} we get an 
$\mrm{L}_{\infty}$ quasi-iso\-mor\-phism
\[ \Gamma(X, \Psi_{\bsym{\sigma}}) = 
\{ \Gamma(X, \Psi_{\bsym{\sigma}; j}) \}_{j \geq 1} :
\Gamma \big( X, \opn{Mix}_{\bsym{U}} 
(\mcal{T}^{}_{\mrm{poly}, X}) \bigl) \to
\Gamma \big( X, \opn{Mix}_{\bsym{U}}
(\mcal{D}^{}_{\mrm{poly}, X}) \bigl) . \]
\end{cor}

\begin{proof}
Theorem \ref{thm5.1} tells us that $\Psi_{\bsym{\sigma}; 1}$ is a
quasi-isomorphisms of complexes of sheaves. By \cite[Theorem 6.2]{Ye1} it
follows that  
\[ \Gamma(X, \Psi_{\bsym{\sigma}; 1}) :
\Gamma \big( X, \opn{Mix}_{\bsym{U}} 
(\mcal{T}^{}_{\mrm{poly}, X}) \bigl) \to
\Gamma \big( X, \opn{Mix}_{\bsym{U}}
(\mcal{D}^{}_{\mrm{poly}, X}) \bigl) \]
is a quasi-isomorphism. 
\end{proof}

\begin{cor} \label{cor5.5}
The data $(\bsym{U}, \bsym{\sigma})$ induces a bijection
\[ \begin{aligned}
\opn{MC}(\Psi_{\bsym{\sigma}}) & :\
\opn{MC} \Bigl( \Gamma \bigl( X, \opn{Mix}_{\bsym{U}}
(\mcal{T}^{}_{\mrm{poly}, X}) \bigl)[[\hbar]]^+ \Bigr) \\
& \qquad \iso \
\opn{MC} \Bigl( \Gamma \bigl( X, \opn{Mix}_{\bsym{U}}
(\mcal{D}^{}_{\mrm{poly}, X}) \bigl)[[\hbar]]^+ \Bigr) .
\end{aligned} \]
\end{cor}

\begin{proof}
Use Corollary \ref{cor11.1} and \cite[Corollary 3.10]{Ye1}.
\end{proof}

Recall that
$\mcal{T}_{\mrm{poly}}(X) = 
\Gamma(X, \mcal{T}_{\mrm{poly}, X})$
and
$\mcal{D}^{\mrm{nor}}_{\mrm{poly}}(X) = 
\Gamma(X, \mcal{D}^{\mrm{nor}}_{\mrm{poly}, X})$; 
and the latter is the DG Lie algebra of global poly differential 
operators that vanish if one of their arguments is $1$.

Suppose $f : X' \to X$ is an \'etale morphism. According to 
\cite[Prposition 4.6]{Ye2} there are DG Lie algebra homomorphisms
$f^* : \mcal{T}^{}_{\mrm{poly}}(X) \to 
\mcal{T}^{}_{\mrm{poly}}(X')$
and
$f^* : \mcal{D}^{\mrm{nor}}_{\mrm{poly}}(X) \to 
\mcal{D}^{\mrm{nor}}_{\mrm{poly}}(X')$. 
These homomorphisms extend to formal coefficients, and we get functions
\[ \mrm{MC} (f^*) : 
\mrm{MC} \bigl( \mcal{T}_{\mrm{poly}}(X)[[\hbar]]^+ \bigr) \to
\mrm{MC} \bigl( \mcal{T}_{\mrm{poly}}(X')[[\hbar]]^+ \bigr) \]
etc.

One says that $X$ is a {\em $\mcal{D}$-affine variety} if 
$\mrm{H}^q(X, \mcal{M}) = 0$ for every quasi-coherent left
$\mcal{D}_X$-module $\mcal{M}$ and every $q > 0$.

\begin{thm} \label{thm5.2}
Let $X$ be an irreducible smooth separated $\K$-scheme. Assume $X$ is 
$\mcal{D}$-affine. Then there is a canonical function 
\[ Q : 
\mrm{MC} \bigl( \mcal{T}_{\mrm{poly}}(X)[[\hbar]]^+ \bigr) 
\to
\mrm{MC} \bigl( \mcal{D}^{\mrm{nor}}_{\mrm{poly}}(X)[[\hbar]]^+ 
\bigr) \]
called the {\em quantization map}. It has the following properties:
\begin{enumerate}
\rmitem{i} The function $Q$ preserves first order terms.
\rmitem{ii} The function $Q$ respects \'etale morphisms. Namely if $X'$ is
another $\mcal{D}$-affine scheme, with quantization map $Q'$,
and if $f : X' \to X$ is an \'etale morphism, then 
\[ Q' \circ \mrm{MC}(f^*) = \mrm{MC}(f^*)  \circ Q . \]
\rmitem{iii} If $X$ is affine, then $Q$ is bijective. 
\rmitem{iv} The function $Q$ 
is characterized as follows. Choose an open covering
$\bsym{U} = \{ U_{0}, \ldots, U_{m} \}$ 
of $X$ consisting of affine open sets, each admitting 
an \'etale coordinate system. Let 
$\bsym{\sigma}$ be the associated simplicial section of the bundle
$\opn{LCC} X \to X$. Then there is a commutative 
diagram 
\[ 
\begin{CD}
\mrm{MC} \bigl( \mcal{T}^{}_{\mrm{poly}}(X)[[\hbar]]^+ \bigr) 
@> Q >>
\mrm{MC} \bigl( \mcal{D}^{\mrm{nor}}_{\mrm{poly}}(X)[[\hbar]]^+ 
\bigr) \\
@V \opn{MC}(\eta_{\mcal{T}}) VV  @V \opn{MC}(\eta_{\mcal{D}}) VV \\
\mrm{MC} \bigl( \Gamma(X, \opn{Mix}_{\bsym{U}}  
(\mcal{T}^{}_{\mrm{poly}, X}))[[\hbar]]^+ \bigr)
@> \opn{MC}(\Psi_{\bsym{\sigma}}) >>
\mrm{MC} \bigl( 
\Gamma(X, \opn{Mix}_{\bsym{U}}
(\mcal{D}^{}_{\mrm{poly}, X}))[[\hbar]]^+ \bigr) 
\end{CD} 
\]
in which the arrows $\opn{MC}(\Psi_{\bsym{\sigma}})$ and
$\opn{MC}(\eta_{\mcal{D}})$ are bijections. 
Here $\Psi_{\bsym{\sigma}}$ is the $\mrm{L}_{\infty}$ 
quasi-isomorphism from Theorem \tup{\ref{thm5.1}}, and 
$\eta_{\mcal{T}}, \eta_{\mcal{D}}$ are the inclusions of DG Lie
algebras.
\end{enumerate}
\end{thm}

Let's elaborate a bit on the statement above. It says that 
to any MC solution
$\alpha = \sum_{j = 1}^{\infty} \alpha_j \hbar^j
\in \mcal{T}^{1}_{\mrm{poly}}(X)[[\hbar]]^+$
there corresponds an MC solution 
$\beta = \sum_{j = 1}^{\infty} \beta_j \hbar^j 
\in \mcal{D}^{\mrm{nor}, 1}_{\mrm{poly}}(X)[[\hbar]]^+$. 
The element $\beta = Q(\alpha)$ 
is uniquely determined up to gauge equivalence by the group
$\opn{exp}(\mcal{D}^{\mrm{nor}, 0}_{\mrm{poly}}(X)[[\hbar]]^+)$.
Given any local sections $f, g \in \mcal{O}_{X}$ one has
\begin{equation} \label{eqn12.9}
\smfrac{1}{2} (\beta_1(f, g) - \beta_1(g, f)) = 
\alpha_1(f, g) \in \mcal{O}_{X} . 
\end{equation}
The quantization map $Q_{}$ can be calculated (at 
least in theory)
using the collection of sections $\bsym{\sigma}$ and the universal 
formulas for deformation in \cite[Theorem 3.13]{Ye1}. 

We'll need a lemma before proving the theorem.

\begin{lem} \label{lem11.2}
Let $f, g \in \mcal{O}_{X} = \mcal{D}^{-1}_{\mrm{poly}, X}$
be local sections.
\begin{enumerate}
\item For any  
$\beta \in 
\opn{Mix}^0_{\bsym{U}}(\mcal{D}^{1}_{\mrm{poly}, X})$
one has 
\[ [[ \beta, f], g] = \beta(g, f) - \beta(f, g) 
\in \opn{Mix}^0_{\bsym{U}}(\mcal{O}_{X}) . \]
\item For any 
$\beta \in 
\opn{Mix}^1_{\bsym{U}}(\mcal{D}^{0}_{\mrm{poly}, X})
\oplus \opn{Mix}^2_{\bsym{U}}(\mcal{D}^{-1}_{\mrm{poly}, X})$
one has $[[ \beta, f], g] = 0$.
\item Let 
$\gamma \in 
\opn{Mix}_{\bsym{U}}(\mcal{D}^{}_{\mrm{poly}, X})^0$,
and define
$\beta := (\d_{\mrm{mix}} + \d_{\mcal{D}})(\gamma)$.
Then $[[ \beta, f], g] = 0$.
\end{enumerate}
\end{lem}

\begin{proof}
(1) \cite[Proposition 6.3]{Ye1} implies that the embedding 
(\cite[(6.1)]{Ye1}:
\[ \begin{aligned}
& \opn{Mix}_{\bsym{U}}(\mcal{D}^{}_{\mrm{poly}, X}) \\
& \qquad \subset \
\bigoplus_{p,q,r} \
{\displaystyle \prod_{j \in \mbb{N}}}\ \
{\displaystyle \prod_{\bsym{i} \in 
\bsym{\Delta}_j^{m}}} \,
g_{\bsym{i} *}\, g_{\bsym{i}}^{-1}\,
\big( \Omega^{q}(\bsym{\Delta}^j_{\K})
\, \what{\otimes} \, (\Omega^{p}_X \otimes_{\mcal{O}_{X}}
\mcal{P}_X \otimes_{\mcal{O}_{X}} 
\mcal{D}^{r}_{\mrm{poly}, X}) \big)
\end{aligned} \]
is a DG Lie algebra homomorphism.
So by continuity we might as well assume that 
$\beta = a D$ with 
$a \in \Omega^0_X = \mcal{O}_X$ and 
$D \in \mcal{D}^{1}_{\mrm{poly}, X}$.
Moreover, since the Lie bracket of 
$\Omega^{}_X \otimes_{\mcal{O}_{X}}
\mcal{P}_X \otimes_{\mcal{O}_{X}} \mcal{D}^{}_{\mrm{poly}, X}$
is $\Omega_X$-bilinear, we may assume that $a = 1$, 
i.e.\ $\beta = D$. Now the assertion
is clear from the definition 
of the Gerstenhaber Lie bracket, see \cite[Section 3.4.2]{Ko}.

\medskip \noindent
(2) Applying the same reduction as above, but with 
$D \in \mcal{D}^{r}_{\mrm{poly}, X}$
and $r \in \set{0, -1}$, we get
$[[D, f], g] \in \mcal{D}^{r-2}_{\mrm{poly}, X} = 0$.

\medskip \noindent
(3) By part (2) it suffices to show that 
$[[ \beta, f], g] = 0$ for
$\beta := \d_{\mcal{D}}(\gamma)$ and
$\gamma \in \opn{Mix}^{0}_{\bsym{U}}
(\mcal{D}^{0}_{\mrm{poly}, X})$.
As explained above we may further assume that 
$\gamma = D \in \mcal{D}^{0}_{\mrm{poly}, X}$. 
Now the formulas for $\d_{\mcal{D}}$ and $[-,-]$
in \cite[Section 3.4.2]{Ko} imply that \linebreak
$[[\d_{\mcal{D}}(D), f], g] = 0$. 
\end{proof}

\begin{proof}[Proof of Theorem \tup{\ref{thm5.2}}]
Step 1. Take an open covering $\bsym{U}$ as in property (iv). 
Since the sheaves 
$\mcal{D}^{\mrm{nor}, p}_{\mrm{poly}, X}$
are quasi-coherent left $\mcal{D}_X$-modules, it follows that
$\mrm{H}^q(X, \mcal{D}^{\mrm{nor}, p}_{\mrm{poly}, X}) = 0$
for all $p$ and all $q>0$.
Therefore 
$\Gamma(X, \mcal{D}^{\mrm{nor}}_{\mrm{poly}, X}) = 
\mrm{R} \Gamma(X, \mcal{D}^{\mrm{nor}}_{\mrm{poly}, X})$
in the derived category $\msf{D}(\cat{Mod} \K)$.
Now by \cite[Theorem 3.12]{Ye1} the inclusion 
$\mcal{D}^{\mrm{nor}}_{\mrm{poly}, X} \to
\mcal{D}^{}_{\mrm{poly}, X}$
is a quasi-isomorphism, and by \cite[Theorem 6.2(1)]{Ye1}
the inclusion 
$\mcal{D}^{}_{\mrm{poly}, X} \to 
\opn{Mix}_{\bsym{U}}(\mcal{D}^{}_{\mrm{poly}, X})$
is a quasi-isomorphism. According to \cite[Theorem 6.2(2)]{Ye1} we have
$\Gamma \bigl( X, \opn{Mix}_{\bsym{U}}
(\mcal{D}^{}_{\mrm{poly}, X}) \bigr) = 
\mrm{R} \Gamma \bigl( X, \opn{Mix}_{\bsym{U}}
(\mcal{D}^{}_{\mrm{poly}, X}) \bigr)$.
The conclusion is that
\begin{equation} \label{eqn7.4}
\mcal{D}^{\mrm{nor}}_{\mrm{poly}}(X) = 
\Gamma(X, \mcal{D}^{\mrm{nor}}_{\mrm{poly}, X}) \to
\Gamma \bigl( X, \opn{Mix}_{\bsym{U}}
(\mcal{D}^{}_{\mrm{poly}, X}) \bigr)
\end{equation}
is a quasi-isomorphism of complexes of $\K$-modules. 
But in view of \cite[Proposition 6.3]{Ye1},
this is in fact a quasi-isomorphism of DG Lie algebras. 

From (\ref{eqn7.4}) we deduce that 
\[ \eta_{\mcal{D}} :
\mcal{D}^{\mrm{nor}}_{\mrm{poly}}(X)[[\hbar]]^+ \to
\Gamma \bigl( X, \opn{Mix}_{\bsym{U}}
(\mcal{D}^{}_{\mrm{poly}, X}) \bigr)[[\hbar]]^+ \]
is a quasi-isomorphism of DG Lie algebras. Using 
\cite[Corollary 3.10]{Ye1} we see that 
$\opn{MC}(\eta_{\mcal{D}})$ is bijective.
Therefore the diagram in property (iv) defines $Q_{}$ uniquely.

According to Corollary \ref{cor5.5}, the arrow marked
$\opn{MC}(\Psi_{\bsym{\sigma}})$ is a bijection. So we have estanlished
property (iv), except for the independence of the open covering. 

\medskip \noindent
Step 2.  The left vertical
arrow comes from the DG Lie algebra homomorphism
\[ \eta_{\mcal{T}} : \mcal{T}_{\mrm{poly}}(X)[[\hbar]]^+ \to
\Gamma \bigl( X, \opn{Mix}_{\bsym{U}}
(\mcal{T}^{}_{\mrm{poly}, X}) \bigr)[[\hbar]]^+ , \]
which is a quasi-isomorphism when 
$\mrm{H}^q(X, \mcal{T}^{p}_{\mrm{poly}, X}) = 0$
for all $p$ and all $q>0$. So in case $X$ is affine, the quantization map $Q$ is
bijective. This establishes property (iii).

\medskip \noindent
Step 3.
Now suppose $\bsym{U}' = \{ U'_{0}, \ldots, U'_{m'} \}$ 
is another such affine open covering of $X$, with sections
$\sigma'_{i} : U'_{i} \to \opn{LCC} X$. 
Without loss of generality we may assume that $m' \geq m$, and 
that $U'_{i} = U_{i}$ and $\sigma'_{i} = \sigma_{i}$
for all $i \leq m$. There is a morphism of simplicial schemes
$f : \bsym{U} \to \bsym{U}'$, that is an open and closed embedding.
Correspondingly there is a commutative diagram 
\[ \begin{CD}
\mrm{MC} \bigl( \mcal{T}^{}_{\mrm{poly}}(X)[[\hbar]]^+ \bigr) 
@> Q >>
\mrm{MC} \bigl( \mcal{D}^{\mrm{nor}}_{\mrm{poly}}(X)[[\hbar]]^+ 
\bigr) \\
@V \opn{MC}(\eta_{\mcal{T}}) VV  @V \opn{MC}(\eta_{\mcal{D}}) VV \\
\mrm{MC} \bigl( 
\Gamma(X, \opn{Mix}_{\bsym{U}'}  
(\mcal{T}^{}_{\mrm{poly}, X}))[[\hbar]]^+ \bigr)
@> \opn{MC}(\Psi_{\bsym{\sigma}'}) >>
\mrm{MC} \bigl( 
\Gamma(X, \opn{Mix}_{\bsym{U}'}
(\mcal{D}^{}_{\mrm{poly}, X}))[[\hbar]]^+ \bigr) \\
@V \mrm{MC}(f^*) VV  @V \opn{MC}(f^*) VV \\
\mrm{MC} \bigl( 
\Gamma(X, \opn{Mix}_{\bsym{U}}
(\mcal{T}^{}_{\mrm{poly}, X}))[[\hbar]]^+ \bigr)
@> \opn{MC}(\Psi_{\bsym{\sigma}}) >>
\mrm{MC} \bigl( 
\Gamma(X, \opn{Mix}_{\bsym{U}}
(\mcal{D}^{}_{\mrm{poly}, X}))[[\hbar]]^+ \bigr)  ,
\end{CD} \]
where the vertical arrows on the right
are bijections. We conclude that $Q$ is 
independent of $\bsym{U}$ and $\bsym{\sigma}$. This concludes the proof of 
property (iv).

\medskip \noindent
Step 4.
Suppose $f : X' \to X$ is an \'etale morphism. Then we can choose an affine
open covering $\bsym{U}'$ of $X'$ that refines $\bsym{U}$ in the obvious
sense. Each of the open sets $U'_i$ inherits an \'etale coordinate system, and
hence a section $\sigma'_i : U'_i \to \opn{LCC} X'$.
We get a commutative diagram 
\[ \begin{CD}
\mrm{MC} \bigl( 
\Gamma(X, \opn{Mix}_{\bsym{U}}  
(\mcal{T}^{}_{\mrm{poly}, X}))[[\hbar]]^+ \bigr)
@> \opn{MC}(\Psi_{\bsym{\sigma}}) >>
\mrm{MC} \bigl( 
\Gamma(X, \opn{Mix}_{\bsym{U}}
(\mcal{D}^{}_{\mrm{poly}, X}))[[\hbar]]^+ \bigr) \\
@V \mrm{MC}(f^*) VV  @V \opn{MC}(f^*) VV \\
\mrm{MC} \bigl( 
\Gamma(X', \opn{Mix}_{\bsym{U}'}
(\mcal{T}^{}_{\mrm{poly}, X'}))[[\hbar]]^+ \bigr)
@> \opn{MC}(\Psi_{\bsym{\sigma}'}) >>
\mrm{MC} \bigl( 
\Gamma(X', \opn{Mix}_{\bsym{U}'}
(\mcal{D}^{}_{\mrm{poly}, X'}))[[\hbar]]^+ \bigr)  ,
\end{CD} \]
This proves property (ii). 

\medskip \noindent
Step 5.
Finally we must show that $Q$ preserves first order terms, i.e.\ property
(i). Let
\[ \alpha = \sum_{j=1}^{\infty} \alpha_j \hbar^j \in 
\mcal{T}_{\mrm{poly}}(X)^1[[\hbar]]^+ \]
be an MC solution, and let 
\[ \beta = \sum_{j=1}^{\infty} \beta_j \hbar^j
\in \mcal{D}^{\mrm{nor}}_{\mrm{poly}}(X)^1[[\hbar]]^+ \]
be an MC solution such that $\beta = Q(\alpha)$ modulo gauge equivalence.
This means that there exists some 
\[ \gamma = \sum_{k \geq 1} \gamma_k \hbar^k \in
\Gamma(X, \opn{Mix}_{\bsym{U}}
(\mcal{D}^{}_{\mrm{poly}, X}))^0[[\hbar]]^+ \]
such that
\[ \sum_{j \geq 1} \smfrac{1}{j!}
\Psi_{\bsym{\sigma}; j}(\alpha^j) = 
\opn{exp}(\opn{af})(\opn{exp}(\gamma))(\beta) , \]
with notation as in \cite[Lemma 3.2]{Ye1}.
Cf.\ \cite[Theorem 3.8]{Ye1}.
In the first order term (i.e.\ the coefficient of $\hbar^1$) of 
this equation we have
\begin{equation} \label{eqn10.5}
\Psi_{\bsym{\sigma}; 1}(\alpha_1) = \beta_1 - 
(\d_{\mrm{mix}} + \d_{\mcal{D}})(\gamma_1) ;
\end{equation}
see \cite[equation (3.3)]{Ye1}. 

In order to apply Lemma \ref{lem11.2}(2), we are interested in in the component
of $\Psi_{\bsym{\sigma}; 1}(\alpha_1)$
living in the summand
$\opn{Mix}^0_{\bsym{U}}(\mcal{D}^{1}_{\mrm{poly}, X})$.
But this is exactly 
\[ \opn{gr}^0_{\mrm{G}}(\Psi_{\bsym{\sigma}; 1})(\alpha_1)
\in \opn{gr}^0_{\mrm{G}} \opn{Mix}_{\bsym{U}}(\mcal{D}^{1}_{\mrm{poly}, X})
= \opn{Mix}^0_{\bsym{U}}(\mcal{D}^{1}_{\mrm{poly}, X}) . \]
Since according to Theorem \ref{thm5.1} we have
\[ \opn{gr}^0_{\mrm{G}}(\Psi_{\bsym{\sigma}; 1}) = 
\opn{gr}^0_{\mrm{G}}(\opn{Mix}^{}_{\bsym{U}}(\mcal{U}_{1})) , \]
it follows that the component we are interested in is
\[ \opn{gr}^0_{\mrm{G}}(\opn{Mix}^{}_{\bsym{U}}(\mcal{U}_{1}))(\alpha_1) = 
\mcal{U}_{1}(\alpha_1) . \]

Now take any two local sections $f, g \in \mcal{O}_{X}$. Using Lemma
\ref{lem11.2} we get
\[ [[\Psi_{\bsym{\sigma}; 1}(\alpha_1), f], g] = 
[[ \mcal{U}_{1}(\alpha_1), f], g] =
\mcal{U}_{1}(\alpha_1)(g, f) - \mcal{U}_{1}(\alpha_1)(f, g) = 
- 2 \alpha_1(f, g) ,  \]
\[ [[ \beta_1, f], g] = \beta_1(g, f) - \beta_1(f, g) \]
and
\[ [[ (\d_{\mrm{mix}} + \d_{\mcal{D}})(\gamma_1), f], g] = 0 . \]
Combining these equations with equation (\ref{eqn10.5}) we see that equation
(\ref{eqn12.9}) indeed holds. So the proof is done.
\end{proof}

\begin{cor} \label{cor9.1}
Let $X$ be an irreducible smooth separated $\K$-scheme. Assume $X$ is 
$\mcal{D}$-affine. Then the quantization map $Q$
of Theorem \tup{\ref{thm5.2}} may be interpreted as a 
canonical function
\[ Q : 
\frac{ \{ \tup{formal Poisson structures on $X$} \} }
{\tup{gauge equivalence}} \to
\frac{ \{ \tup{deformation quantizations of $\mcal{O}_X$} \} }
{\tup{gauge equivalence}} . \]
The quantization map $Q$ preserves first order terms, and commutes with
\'etale morphisms 
$f : X' \to X$. If $X$ is affine then $Q$ is bijective.
\end{cor}

\begin{proof}
By definition the left side is 
$\mrm{MC} \bigl( \mcal{T}^{}_{\mrm{poly}}(X)[[\hbar]]^+ \bigr)$.
On the other hand, according to \cite[Theorem 1.13]{Ye1} every 
deformation quantization of $\mcal{O}_{X}$ can be trivialized 
globally, and by \cite[Proposition 1.14]{Ye1} any gauge equivalence 
between globally trivialized deformation quantizations is a
global gauge equivalence. Hence the right side is 
$\mrm{MC} \bigl( \mcal{D}^{\mrm{nor}}_{\mrm{poly}}(X)[[\hbar]]^+ 
\bigr)$.
\end{proof}

\section{Miscellaneous Errors}

Here is a list of minor errors in the paper \cite{Ye1}. 

\begin{enumerate}
\item Section 3, bottom of page 395: the formula should be
\[ \opn{af}(\gamma)(\omega) := [\gamma, \omega] - \mrm{d}(\gamma) 
 = \opn{ad}(\gamma)(\omega) - \d(\gamma) \in
\mfrak{m} \hatotimes{} \mfrak{g}^1  , \]
\item Definition 5.2, page 411: the formula should be 
\[ \nabla_{\mcal{P}} : 
\mcal{P}_X \to \Omega^1_X \otimes_{\mcal{O}_{X}}
\mcal{P}_X  . \]

\end{enumerate}

\end{document}